# TRIPLETS OF TRI-HOMOLOGICAL TRIANGLES

Prof. Ion Pătrașcu, The "Frații Buzești" National College, Craiova, Romania
Prof. Florentin Smarandache, The University of New Mexico, Gallup, USA

In this article will prove some theorems in relation to the triplets of homological triangles two by two. These theorems will be used later to build triplets of triangles two by two tri-homological.

**I     Theorems on the triplets of homological triangles**
**Theorem 1**

Two triangles are homological two by two and have a common homological center (their homological centers coincide) then their homological axes are concurrent.

**Proof**

Let's consider the homological triangles $A_1B_1C_1$, $A_2B_2C_2$, $A_3B_3C_3$ whose common homological center is $O$ (see figure 1.)

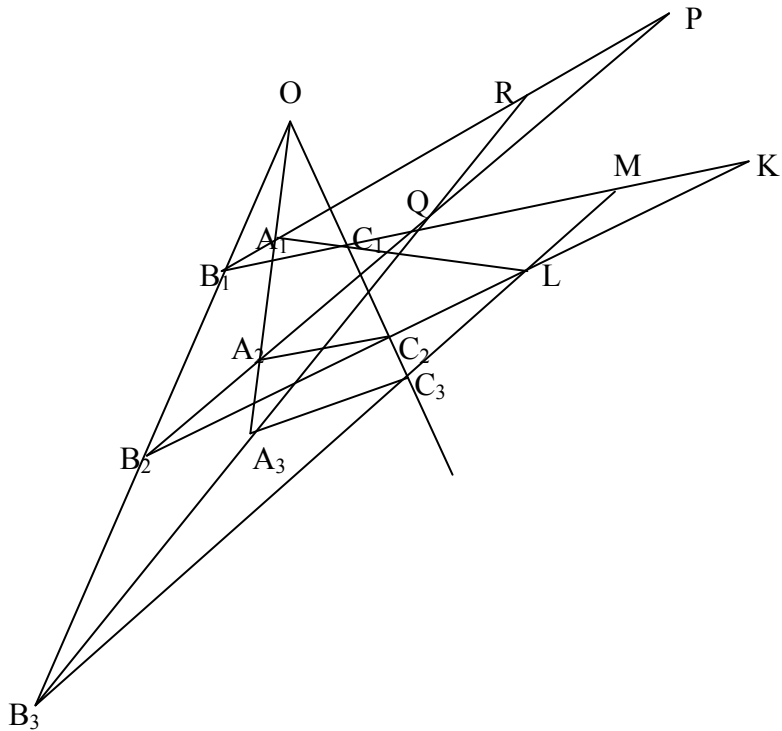

Fig. 1

We consider the triangle formed by the intersections of the lines: $A_1B_1$, $A_2B_2$, $A_3B_3$ and we note it $PQR$ and the triangle formed by the intersection of the lines $B_1C_1$, $B_2C_2$, $B_3C_3$ and we'll note it $KLM$.

We observe that $PR \cap KM = \{B_1\}, RQ \cap ML = \{B_2\}$, $PQ \cap KL=\{B_3\}$ and because $B_1, B_2, B_3$ are collinear it results, according to the Desargues reciprocal theorem that the triangles $PQR$ and $KLM$ are homological, therefore $PK, RM, QL$ are concurrent lines.

The line $PK$ is the homological axes of triangles $A_1B_1C_1$ and $A_2B_2C_2$, the line $RM$ is the homological axis for triangles $A_1B_1C_1$ and $A_3B_3C_3$, and the line $QL$ is the homological axis for triangles $A_2B_2C_2$ and $A_3B_3C_3$, which proves the theorem.

**Remark 1**

Another proof of this theorem can be done using the spatial vision; if we imagine figure 1 as being the correspondent of a spatial figure, we notice that the planes $(A_1B_1C_1)$ and $(A_2B_2C_2)$ have in common the line $PK$, similarly the planes $(A_1B_1C_1)$ and $(A_3B_3C_3)$ have in common the line $QL$. If $\{O'\} = PK \cap LQ$ then $O'$ will be in the plane $(A_2B_2C_2)$ and in the plane $(A_3B_3C_3)$, but these planes intersect by the line $RM$, therefore $O'$ belongs to this line as well. The lines $PK, RM, QL$ are the homological axes of the given triangles and therefore these are concurrent in $O'$.

**Theorem 2**

If three triangles are homological two by two and have the same homological axis (their homological axes coincide) then their homological axes are collinear.

**Proof**

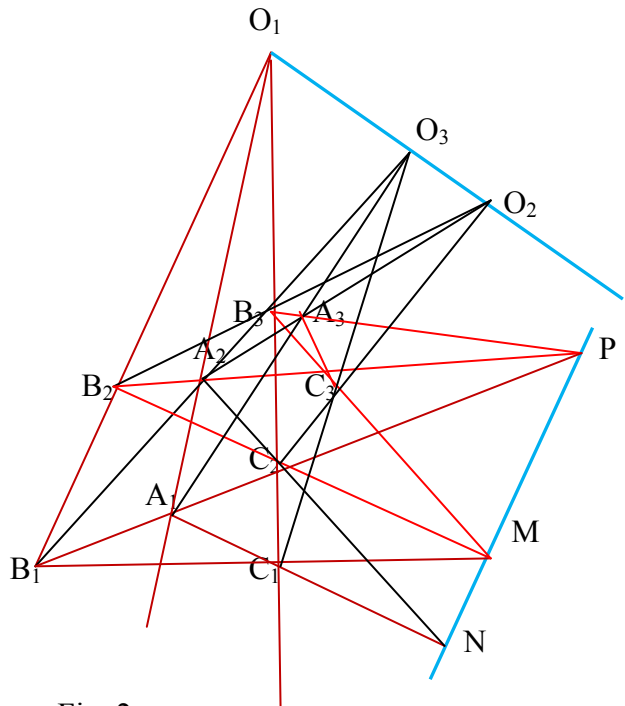

Fig. 2

Let's consider the homological triangles two by two $A_1B_1C_1$, $A_2B_2C_2$, $A_3B_3C_3$. We note $M, N, P$ their common homological axis (see figure 2). We note $O_1$ the homological center of the triangles $A_1B_1C_1$ and $A_2B_2C_2$, with $O_2$ the homological center of the triangles $A_2B_2C_2$ and $A_3B_3C_3$ and with $O_3$ the homological center of the triangles $A_3B_3C_3$ and $A_1B_1C_1$.

We consider the triangles $A_1A_2A_3$ and $B_1B_2B_3$, and we observe that these are homological because $A_1B_1$, $A_2B_2$, $A_3B_3$ intersect in the point $P$ which is their homological center. The homological axis of these triangles is determined by the points
, $\{O_1\} = A_1A_2 \cap B_1B_2$, $\{O_2\} = A_2A_3 \cap B_2B_3$, $\{O_3\} = A_1A_3 \cap B_1B_3$
therefore the points $O_1, O_2, O_3$ are collinear and this concludes the proof of this theorem.

**Theorem 3** (The reciprocal of theorem 2)
If three triangles are homological two by two and have their homological centers collinear, then these have the same homological axis.

**Proof**
We will use the triangles from figure 2. Let therefore $O_1, O_2, O_3$ the three homological collinear points. We consider the triangles $B_1B_2B_3$ and $C_1C_2C_3$, we observe that these admit as homological axis the line $O_1O_2O_3$.
Because
$$\{O_1\} = B_1B_2 \cap C_1C_2, \{O_2\} = B_2B_3 \cap C_2C_3, \{O_3\} = B_1B_3 \cap C_1C_3,$$

It results that these have as homological center the point $\{M\} = B_1C_1 \cap B_2C_2 \cap B_3C_3$.

Similarly for the triangles $A_1A_2A_3$ and $C_1C_2C_3$ have as homological axis $O_1O_2O_3$ and the homological center $M$. We also observe that the triangles $A_1A_2A_3$ and $B_1B_2B_3$ are homological and $O_1O_2O_3$ is their homological axis, and their homological center is the point $P$. Applying the theorem 2, it results that the points $M, N, P$ are collinear, and the reciprocal theorem is then proved.

**Theorem 4** (The Veronese theorem)
If the triangles $A_1B_1C_1$, $A_2B_2C_2$ are homological and
$\{A_3\} = B_1C_2 \cap B_2C_1, \{B_3\} = A_1C_2 \cap A_2C_1, \{C_3\} = A_1B_2 \cap A_2B_1$
then the triangle $A_3B_3C_3$ is homological with each of the triangles $A_1B_1C_1$ and $A_2B_2C_2$, and their homological centers are collinear.

**Proof**
Let $O_1$ be the homological center of triangles $A_1B_1C_1$ and $A_2B_2C_2$ (see figure 3) and $A', B', C'$ their homological axis.

We observe that $O_1$ is a homological center also for the triangles $A_1B_1C_2$ and $A_2B_2C_1$. The homological axis of these triangles is $C', A_3, B_3$. Also $O_1$ is the homological center for the triangles

$B_1C_1A_2$ and $B_2C_2A_1$, it results that their homological axis is $A', B_3, C_3$

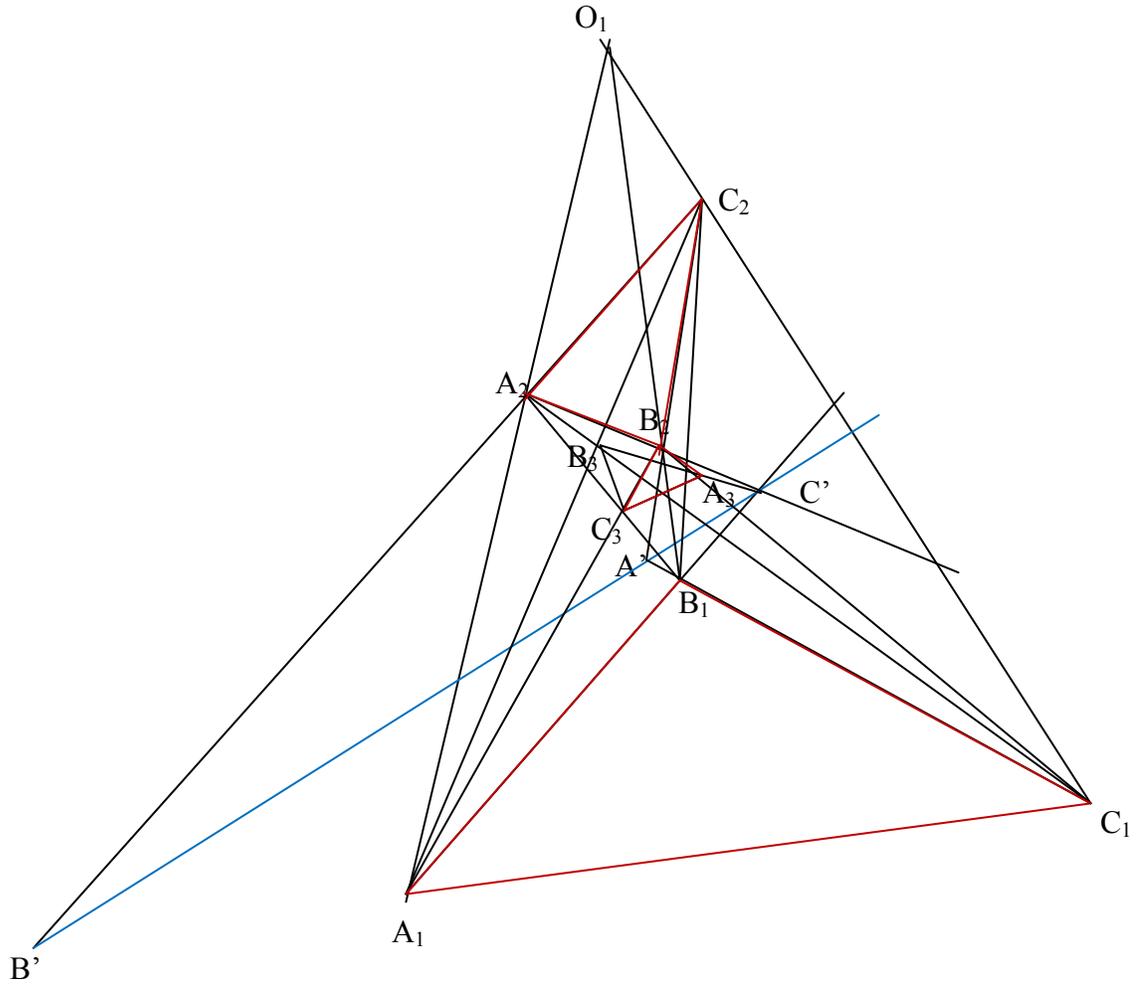

Fig. 3

Similarly, we obtain that the points $B', A_3, C_3$ are collinear, these being on a homological axis of triangle $C_1A_1B_2$ and $C_2A_2B_1$. The triplets of the collinear points $(C', A_3, B_3)$, $(B', A_3, C_3)$ and $(A', B_3, C_3)$ show that the triangle $A_3B_3C_3$ is homological with triangle $A_1B_1C_1$ and with the triangle $A_2B_2C_2$.

The triangles $A_1B_1C_1$, $A_2B_2C_2$, $A_3B_3C_3$ are homological two by two and have the same homological axis $A', B', C'$. Using theorem 3, it results that their homological centers are collinear points.

## II. Double-homological triangles
### Definition 1

We say that the triangles $A_1B_1C_1$ and $A_2B_2C_2$ are double-homological or bi-homological if these are homological in two modes.

**Theorem 5**

Let's consider the triangles $A_1B_1C_1$ and $A_2B_2C_2$ such that

$$B_1C_1 \cap B_2C_2 = \{P_1\}, B_1C_1 \cap A_2C_2 = \{Q_1\}, B_1C_1 \cap A_2B_2 = \{R_1\}$$

$$A_1C_1 \cap A_2C_2 = \{P_2\}, A_1C_1 \cap A_2B_2 = \{Q_2\}, A_1C_1 \cap B_2C_2 = \{R_2\}$$

$$A_1B_1 \cap A_2B_2 = \{P_3\}, A_1B_1 \cap B_2C_2 = \{Q_3\}, A_1B_1 \cap C_2A_2 = \{R_3\}$$

Then:

$$\frac{P_1B_1 \cdot P_2C_1 \cdot P_3A_1}{P_1C_1 \cdot P_2A_1 \cdot P_3B_1} \cdot \frac{Q_1B_1 \cdot Q_2C_1 \cdot Q_3A_1}{Q_1C_1 \cdot Q_2A_1 \cdot Q_3B_1} \cdot \frac{R_1B_1 \cdot R_2C_1 \cdot R_3A_1}{R_1C_1 \cdot R_2A_1 \cdot R_3B_1} = 1 \qquad (1)$$

**Proof**

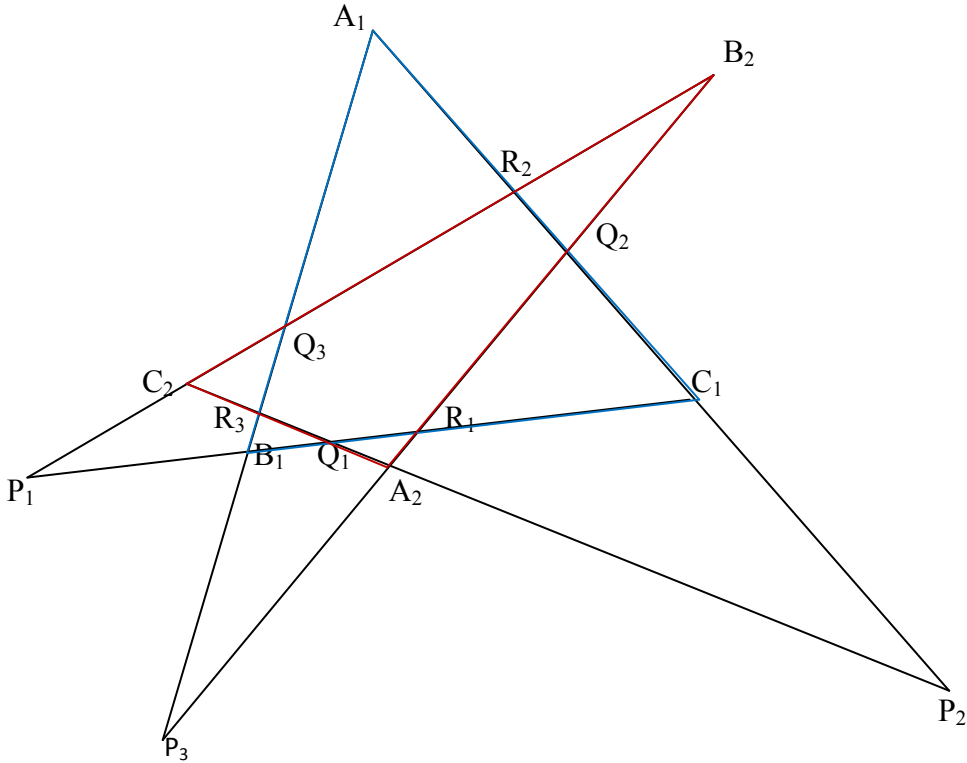

Fig. 4

We'll apply the Menelaus' theorem in the triangle $A_1B_1C_1$ for the transversals $P_1Q_3R_2$, $P_2Q_1R_3$, $P_3Q_2R_1$, (see figure 4).

We obtain

$$\frac{P_1B_1 \cdot R_2C_1 \cdot Q_3A_1}{P_1C_1 \cdot R_2A_1 \cdot Q_3B_1} = 1$$

$$\frac{P_2C_1 \cdot Q_1B_1 \cdot R_3A_1}{P_2A_1 \cdot Q_1C_1 \cdot R_3B_1} = 1$$

$$\frac{P_3A_1 \cdot R_1B_1 \cdot Q_2C_1}{P_3B_1 \cdot R_1C_1 \cdot Q_2A_1} = 1$$

Multiplying these relations side by side and re-arranging the factors, we obtain relation (1).

**Theorem 6**

The triangles $A_1B_1C_1$ and $A_2B_2C_2$ are homological (the lines $A_1A_2$, $B_1B_2$, $C_1C_2$ are concurrent) if and only if:

$$\frac{Q_1B_1 \cdot Q_2C_1 \cdot Q_3A_1}{Q_1C_1 \cdot Q_2A_1 \cdot Q_3B_1} = \frac{R_1C_1 \cdot R_2A_1 \cdot R_3B_1}{R_1B_1 \cdot R_2C_1 \cdot R_3A_1} \qquad (2)$$

**Proof**

Indeed, if $A_1A_2$, $B_1B_2$, $C_1C_2$ are concurrent then the points $P_1, P_2, P_3$ are collinear and the Menelaus' theorem for the transversal $P_1P_2P_3$ in the triangle $A_1B_1C_1$ gives:

$$\frac{P_1B_1 \cdot P_2C_1 \cdot P_3A_1}{P_1C_1 \cdot P_2A_1 \cdot P_3B_1} = 1 \qquad (3)$$

This relation substituted in (1) leads to (2)

**Reciprocal**

If the relation (2) takes place then substituting it in the relation (1) we obtain (3) which shows that $P_1, P_2, P_3$ is the homology axis of the triangles $A_1B_1C_1$ and $A_2B_2C_2$.

**Remark 2**

If in relation (1) two fractions are equal to 1, then the third fraction will be equal to 1, and this leads to the following:

**Theorem 7**

If the triangles $A_1B_1C_1$ and $A_2B_2C_2$ are homological in two modes (are double-homological) then these are homological in three modes (are tri-homological).

**Remark 3**

The precedent theorem can be formulated in a different mod that will allow us to construct tri-homological triangles with a given triangle and of some tri-homological triangles.
Here is the theorem that will do this:

**Theorem 8**

(i) Let $ABC$ a given triangle and $P, Q$ two points in its plane such that $BP$ intersects $CQ$ in $A_1$, $CP$ intersects $AQ$ in $B_1$ and $AP$ intersects $BQ$ in $C_1$.

Then $AA_1, BB_1, CC_1$ intersect in a point $R$.

(ii) If $\cap CP = \{A_2\}, CQ \cap AP = \{B_2\}$, $BP \cap AQ = \{C_2\}$ then the triangles $ABC, A_1B_1C_1, A_2B_2C_2$ are two by two homological and their homological centers are collinear.

**Proof**

(i). From the way how we constructed the triangle $A_1B_1C_1$, we observe that $ABC$ and $A_1B_1C_1$ are double homological, their homology centers being two given points $P, Q$ (see figure 5). Using theorem 7 it results that the triangles $ABC$, $A_1B_1C_1$ are tri-homological, therefore $AA_1, BB_1, CC_1$ are concurrent in point noted R.

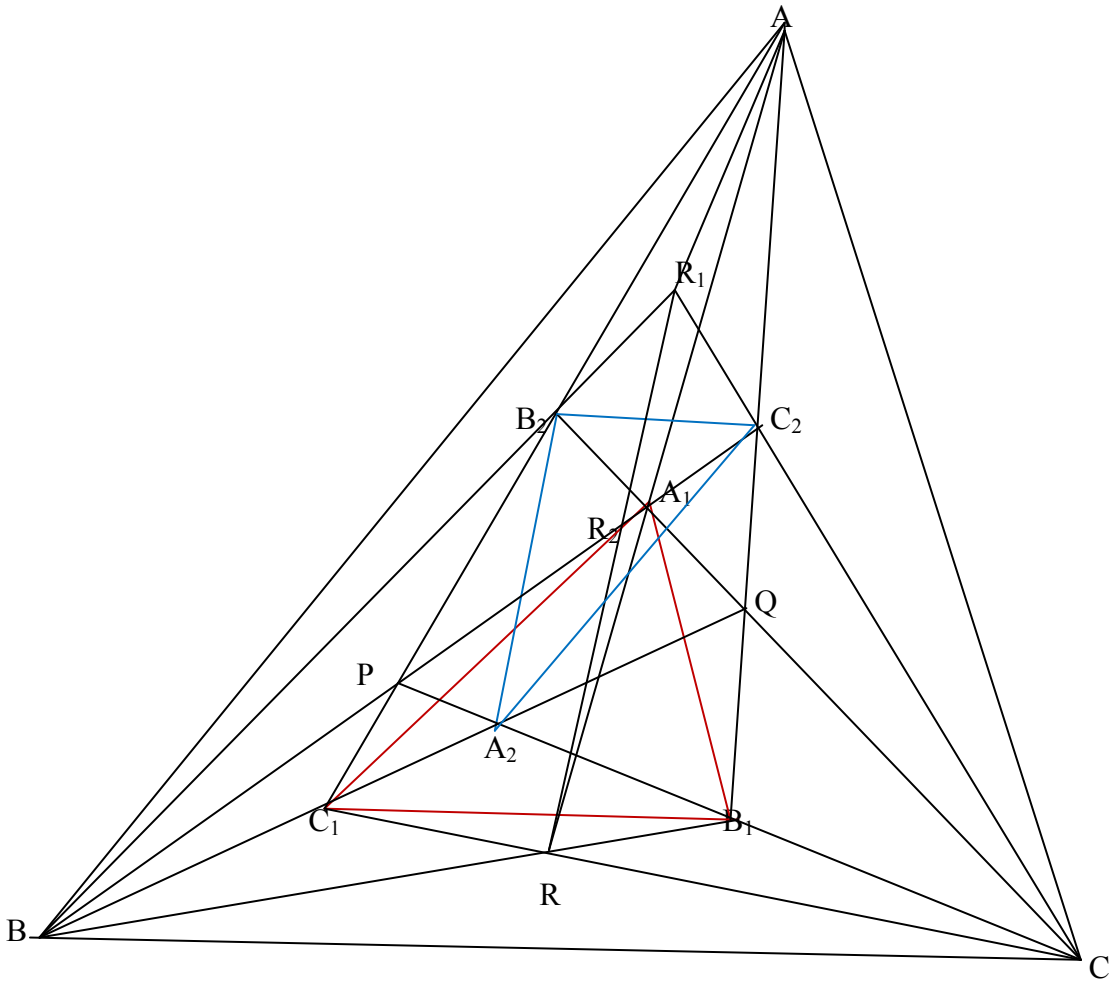

Fig. 5

(ii) The conclusion results by applying the Veronese theorem for the homological triangles $ABC$, $A_1B_1C_1$ that have as homological center the point $R$.

**Remark 4**

We observe that the triangles $ABC$ and $A_2B_2C_2$ are bi-homological, their homological centers being the given points $P, Q$. It results that these are tri-homological and therefore $AA_2, BB_2, CC_2$ are concurrent in the third homological center of these triangle, which we'll note $R_1$.

Similarly we observe that the triangles $A_1B_1C_1$, $A_2B_2C_2$ are double homological with the homological centers $P,Q$; it results that these are tri-homological, therefore $A_1A_2, B_2B_2, C_2C_2$ are concurrent, their concurrence point being notated with $R_2$. In accordance to the Veronese's theorem, applied to any pair of triangles from the triplet $(ABC, A_1B_1C_1, A_2B_2C_2)$ we find that the points $R, R_1, R_2$ are collinear.

**Remark 5**
Considering the points $P, R$ and making the same constructions as in theorem 8 we obtain the triangle $A_3B_3C_3$ which along with the triangles $ABC, A_1B_1C_1$ will form another triplet of triangles tri-homological two by two.

**Remark 6**
The theorem 8 provides us a process of getting a triplet of tri-homological triangles two by two beginning with a given triangle and from two given points in its plane. Therefore if we consider the triangle $ABC$ and as given points the two points of Brocard $\Omega\Omega$ and $\Omega'$, the triangle $A_1B_1C_1$ constructed as in theorem 8 will be the first Brocard's triangle and we'll find that this is a theorem of J. Neuberg: the triangle $ABC$ and the first Brocard triangle are tri-homological. The third homological center of these triangles is noted $\Omega''$ and it is called the Borcard's third point and $\Omega''$ is the isometric conjugate of the simedian center of the triangle $ABC$

**Open problems**
1)   If $T_1, T_2, T_3$ are triangles in a plane, such that $(T_1, T_2)$ are tri-homological, $(T_2, T_3)$ are tri-homological, then are the $(T_1, T_3)$ tri-homological?

2)   If $T_1, T_2, T_3$ are triangles in a plane such that $(T_1, T_2)$ are tri-homological, $(T_2, T_3)$ are tri-homological, $(T_1, T_3)$ are tri-homological and these pairs of triangles have in common two homological centers, then are the three remaining non-common homological centers collinear?